# Bridging probability and calculus: the case of continuous distributions and integrals at the secondary-tertiary transition

Charlotte Derouet[1], Gaetan Planchon[2], Thomas Hausberger[2], and Reinhard Hochmuth[3]

[1]LISEC team AP2E, University of Strasbourg, France, charlotte.derouet@espe.unistra.fr; [2]Institut Montpelliérain Alexander Grothendieck, CNRS, University of Montpellier, France; [3]Leibniz University of Hannover, Germany

*This paper focuses on two mathematical topics, namely continuous probability distributions (CPD) and integral calculus (IC). These two sectors that are linked by the formula $P(a \leq X \leq b) = \int_a^b f(x)dx$ are quite compartmented in teaching classes in France. The main objective is to study whether French students can mobilize the sector of IC to solve tasks in CPD and vice versa at the transition from high school to higher education. Applying the theoretical framework of the Anthropological Theory of the Didactic (ATD), we describe a reference epistemological model (REM) and use it to elaborate a questionnaire in order to test the capacity of students to bridge CPD and IC at the onset of university. The analysis of the data essentially confirms the compartmentalisation of CPD and IC.*

*Keywords: Transition to and across university mathematics; Teaching and learning of analysis and calculus; Teaching and learning of probability; Anthropological Theory of the Didactic.*

## INTRODUCTION

Continuous probability distributions (CPD) and integral calculus (IC) are two topics that are taught in France during the last year of high school (grade 12 of the scientific track). They constitute two sectors (in the sense of the Anthropological Theory of the Didactics, ATD) that belong to the two different but closely related mathematical domains of probability theory and calculus respectively. Indeed, the continuous probability of an event and the definite integral with respect to a non-negative function are both defined as areas of suitable two-dimensional domains in the syllabus (in France), and the formula $P(a \leq X \leq b) = \int_a^b f(x)dx$ is the key for solving several standard tasks in CPD, where $X$ represents a random variable and $f$ its associated density function.

IC is the focus of extended studies in mathematics education: for instance, research (Schneider, 1992; Tran Luong, Bessot & Dorier, 2010; Haddad, 2013) was conducted in the context of Belgian, French, Vietnamese and Tunisian secondary education, including the secondary-tertiary transition (Haddad, 2013). By contrast, there is hardly any literature on CPD and, to make things worse, the available studies mostly put the emphasis on normal distributions (Batanero, Tauber & Meyer, 1999; Wilensky, 1997; Batanero, Tauber & Sánchez, 2004; Pfannkuch & Reading, 2006). Therefore, the teaching-learning phenomena generated by the interrelationship between CPD and IC are still to be investigated.

A first stone was laid by Derouet and Parzysz (2016; see also Derouet, 2016), who studied possible ways to introduce the density function at grade 12 so that students may construct this concept starting from considerations regarding histograms and therefore might relate continuous probability to the integral. By an analysis of textbooks, Derouet could show that the two sectors CPD and IC are very much compartmentalised in the French curriculum (Derouet, 2016, pp. 127-190). For instance, the above formula is seldom justified by a thorough discussion of the definitions involving areas, which certainly hinders the bridging of the two sectors by students.

In this paper, we regard this "compartmentalisation" of knowledge as an institutional phenomenon and therefore use ATD as the theoretical framework (see below). Our goal is to study the impact of this compartmentalisation on the learning of mathematics: are French students able to mobilize the sector of IC to solve tasks in the sector of CPD and vice versa at the transition from high school to higher education?

After the presentation of theoretical constructs from ATD used in this research, we will describe the reference epistemological model (REM) that we elaborated for the types of tasks in CPD and IC with regard to studying interrelations of the two sectors. We will then describe our methodology that builds on the elaboration of a questionnaire, based on the REM, that has been submitted to students at the entrance of university. We finally present results of a primary analysis of the data from the questionnaire and draw some conclusions and perspectives opened up through this study.

**THEORETICAL CONSTRUCTS**

ATD "postulates that any activity related to the production, diffusion or acquisition of knowledge should be interpreted as an ordinary human activity, and thus proposes a general model of human activity built on the key notion of praxeology" (Bosch & Gascon, 2014). The praxeology $\Pi$ is represented by a quadruple $[T/\tau/\theta/\Theta]$: its praxis part (or know-how) consists of a type of tasks $T$ together with a corresponding technique $\tau$ (useful to carry out the tasks $t \in T$ in the scope of $\tau$). The logos part (or know-why) includes two levels of description and justification: the technology $\theta$, i.e. a discourse on the technique, and the theory $\Theta$, which often unifies several technologies.

The elaboration of a reference epistemological model (Florensa, Bosch, & Gascon, 2015) as sequences of praxeologies, for a given body of knowledge, is an important step in any research carried out in the ATD framework. It is the tool that will be used by the researcher to describe, analyse, put in question or design the specific contents that are at the core of a teaching and learning process. In order to build such a model, "mathematical praxeologies are described using data from the different institutions participating in the didactic transposition process, thus including historical, semiotic and sociological research, assuming the institutionalized and socially articulated nature of praxeologies" (loc. cit. p. 2637).

Our study relies on an overview of standard textbooks used at grade 12 in France, as well as the official syllabus, in order to identify the standard praxeologies in CPD and IC that may be related and the nature of this relationship at the praxeological level. An epistemological investigation of the historical development and the interrelation of both domains have previously been carried out in (Derouet, 2016, pp. 67-85). In order to test the effect of the institutional compartmentalisation of knowledge on the learning of mathematics, we need to check the availability of the identified praxeologies in the *praxeological equipment* of students, and then submit tasks to students which need to bridge CPD and IT as mutually interdependent sectors that share techniques or technologies, borrow them from or lend them to each other. Special care must be taken in the phrasing of these bridging tasks, taking into account the effect of *ostensives* (Bosch & Chevallard, 1999), that is to say the role of signs. Indeed, ostensives contribute to the activation of the specific sectors to which they belong and therefore direct students toward specific techniques.

**REFERENCE EPISTEMOLOGICAL MODEL**

Even though problems of quadratures arose in ancient Greece, IC finds its roots as a systematic method in the 17th century. The emergence of continuous probability may be situated in the 18th century with the theory of errors in physical measurements. Various functions were introduced to model the distribution of errors and the area under the curve permitted to evaluate the "theoretical frequency" (so the probability) of the deviation from the "true" value. CPD was thus naturally connected to IC in its historical roots. The gaussian distribution was proposed later by Gauss in 1809.

To identify the different praxeologies, we have analysed 12 textbooks of the grade 12 of the scientific track (edition 2012). We focused on the exercises with a given solution in the textbooks to have access to the usual techniques for the different tasks.

In our study, we will focus on two main types of mathematical tasks $T_I$ and $T_P$, which are related to the mathematical domains of integral calculus and the continuous probability respectively:

- $T_I$: compute a value for an integral $\int_a^b f(x)dx$ for a positive continuous function $f$;
- $T_P$: determine the probability $P(a \leq X \leq b)$ for a random variable $X$ endowed with a density function $f$.

The type of tasks $T_I$ may further be split into two subtypes of tasks, depending on the expected result: an exact value ($T_{I,exact}$) or an approximation ($T_{I,approx}$).

The most useful technique $\tau_{I,exact}$ to solve $T_{I,exact}$ is to compute a primitive of the function and apply the *fundamental theorem of calculus*. The corresponding technology $\theta_{I,exact}$ is given by the fundamental theorem of calculus that relates integrals and primitives: $\int_a^b f(x)dx = F(b) - F(a)$ with $F' = f$. The theory $\Theta_{I,exact}$ includes the definition of the definite integral for a continuous positive function as an area and properties of areas that may be formalised into a local axiomatic theory[1].

The technique $\tau_{I,exact}$ thus resorts to praxeologies dedicated to the computation of primitives. The standard technique at high school level is to use the "tabular of primitives" (deduced from the tabular of derivatives). The technology comprises the properties of the derivative and the theory is that of differential calculus. In the case of piecewise affine functions, an alternative technique to $\tau_{I,exact}$ is to interpret the integral as the area of an elementary surface (or a union of these).

The type of tasks $T_{I,approx}$ may be solved using two main techniques: using a calculator (or software), more or less a blackbox, or applying the "rectangle method". The latter technique $\tau_{I,approx}$ consists in considering the integral as an area, taking a subdivision of the interval of integration and computing the sum of rectangular areas. The technology $\theta_{I,approx}$ comprises the definition of the integral and properties of areas. A further theoretical level $\Theta_{I,approx}$ is mainly non-existent at high school level (cf. endnote 1).

Regarding the type of tasks $T_P$, two cases need to be distinguished, depending on whether a primitive of the density function $f$ is known ($T_{P,prim}$) or not. The latter case is reduced to that of the normal distribution ($T_{P,norm}$), which is dealt with using the implementation of such a distribution in a calculator or software. Computer scientific tools are mainly used as a blackbox by students, which hinders the possibility for students to make connections with IC.

The generic technique $\tau_{P,prim}$ for $T_{P,prim}$ is to compute $\int_a^b f(x)dx$, in other words to resort to the praxeology $\Pi_{I,exact}$. The technology $\theta_{P,prim}$ is given by the formula $P(a \leq X \leq b) = \int_a^b f(x)dx$, and the theory $\Theta_{P,prim}$ comprises the definition of a continuous probability (as the area of the corresponding domain) and the definition of a probability density function. At high school, two particular cases are emphasised and lead to local techniques, as concrete formulas are available for $P(a \leq X \leq b)$ in the case of the exponential and uniform distributions. For instance, the technique $\tau_{P,exp}$ in the case of the exponential may be reduced to computing $e^{\lambda a} - e^{\lambda b}$ with the technological argument $P(a \leq X \leq b) = e^{\lambda a} - e^{\lambda b}$.

Let us recall that our model is based on the study of standard textbooks used in grade 12 classes in France and is dedicated to the description of the teaching-learning of CPD and IC as it actually is (we are not planning task-design at this stage of the research).

In this model, we note the following links between CPD and IC: at the level of the theoretical blocks, praxeologies in both sectors are anchored on the empirical notion of area. At the level of the technique, $\Pi_{P,prim}$ uses $\Pi_{I,exact}$, so that, from an ecological point of view (Bosch & Gascon, 2014, p. 72), CPD contributes to the thriving of such IC praxeologies. By contrast, $\Pi_{I,approx}$ does not seem to be reinvested in CPD (whereas the computation of a probability of the Gaussian distribution could be an opportunity to mobilize $\Pi_{I,approx}$). Conversely, we did not detect elements of the praxis of CPD in the IC sector. This isn't a surprise: IC is regarded as a prerequisite to CPD and precedes the teaching of CPD in all textbooks. Nevertheless, the normal distribution is a prototypical example of a function whose primitive cannot be expressed in terms of available elementary functions. This fact explains the choice of techniques in $\Pi_{P,norm}$ and contributes also to the logos of $\Pi_{I,exact}$ (by complementing the statement that every continuous function admits a primitive). What about the type of tasks $T_{I,norm}$: compute $\int_a^b \frac{e^{-(x-m)^2/2\sigma^2}}{\sigma\sqrt{2\pi}}\,dx$? It could appear in the IC sector through $\Pi_{I,approx}$ only but we didn't find it in any textbook and it is never stated as such in the CPD sector. The most efficient technique requires to use the formula $P(a \leq X \leq b) = \int_a^b f(x)dx$ from right to left: although the equality is symmetrical as an equivalence relation, it isn't symmetrical as a sign which denotes a succession of operations in performing a computation. How would students react to such a task that asks to bridge CPD and IC in an unusual way? This question came to us as a starting point for the elaboration of our questionnaire dedicated to the investigation of the educational effects of the institutional phenomenon of compartmentalisation of knowledge, in the case of CPD and IC.

## THE QUESTIONNAIRE

Our main goal is to test whether students are able or not to connect CPD and IC, and especially mobilize the CPD sector to solve an IC task when ostensives do not indicate explicitly the probability domain. To do so, we have elaborated a questionnaire both to check the availability of standard praxeologies of CPD and IC in the praxeological equipment of students and the capacity of students to complete such bridging tasks.

Bridging tasks appear at the very end of the questionnaire and are stated as follows:

> **Question 6:** Expliquer toutes les méthodes que vous pouvez utiliser pour déterminer une valeur exacte et/ou approchée de l'intégrale suivante : $I = \int_{-0,5}^{1} \frac{e^{-x^2/2}}{\sqrt{2\pi}}\,dx$. On pourra se limiter à donner une idée de la méthode si sa mise en oeuvre est trop compliquée.
>
> **Question 7:** Soit $A$ la fonction définie par $A(\lambda) = \int_0^\lambda f(x)dx$ avec $f(x) = xe^{-x}$, pour tout $\lambda \in [0;+\infty[$. On peut démontrer que $\lim_{\lambda \to +\infty} \int_0^\lambda xe^{-x}dx = 1$. D'après ce résultat, expliciter tout ce que vous pouvez dire sur la fonction $A$ et la fonction $f$.
>
> **Translation:**

> **Question 6:** Explain all the methods that can be used to determine an exact and/or approximate value of the following integral: $I = \int_{-0,5}^{1} \frac{e^{-x^2/2}}{\sqrt{2\pi}} dx$. It suffices to give an idea of the method if its implementation is too complicated.
>
> **Question 7:** Let $A$ be the function defined by $A(\lambda) = \int_0^\lambda f(x)dx$ with $f(x) = xe^{-x}$, for all $\lambda \in [0; +\infty[$. It can be proved that $\lim_{\lambda \to +\infty} \int_0^\lambda xe^{-x}dx = 1$. According to this result, note everything that you can say about the function $A$ and the function $f$.

*Figure 1: bridging tasks submitted to students*

Question 6 contains an instance of the type of tasks $T_{I,norm}$ discussed in the REM. The task is stated in an opened way, asking for every method that students may know to compute an exact or approximate value for the Gaussian integral. Praxeologies $\Pi_{I,exact}$ and $\Pi_{I,approx}$ should therefore also show up. In question 7, we intend to check whether students can say that $A$ is both a primitive for $f$ and a probability associated with the density function $f$, or restrict to the IC sector with an interpretation in terms of areas.

Previous questions intend to "activate" both sectors CPD and IC equally. In this respect, question 1 offers a routine task of type $T_{I,exact}$ in the case of a straightforward exponential function. Analogously, the first part of question 4 is a routine task of type $T_{P,exp}$ (compute $P(1 \leq X \leq 5)$) when $X$ has an exponential distribution of parameter 3). In its second part, students are asked for a graphical interpretation of the probability $P(1 \leq X \leq 5)$ that is in fact defined as an area in the high-school syllabus, as well as the integral: this interpretation is therefore essential in order to link the logos of $\Pi_P$ and that of $\Pi_I$.

Question 2 activates the CPD sector by soliciting an element of the logos of $\Pi_P$, namely the properties that define a density function. This logos is crucial in the bridging question 7, which is stated in the IC sector without any reference to CPD. Question 3 tests if students are able to retrieve the definitions of both the exponential and normal distributions by specifically asking for those in the case of simple parameters (the reduced centered gaussian law). The latter is an element of the logos of $\Pi_{P,norm}$: we wish to check if students are able to recognize the normal distribution in the statement of the bridging question 6, while taking care not to direct them towards a specific technique (hence the order of questions).

In question 5, we rather activate the IC sector, more precisely elements of the logos of $\Pi_{I,exact}$ (primitives), but we have in mind the praxeology $\Pi_{P,norm}$ in relation to the bridging question 6: students are asked to provide an example of a continuous function, if it exists (or justify the impossibility), that a) doesn't possess primitives b) admits primitives but expressions for these are not "explicitly known" (expected answer: the density function of the normal law).

Summarizing, by the questions 1 to 4 we want to investigate whether the students master techniques $\tau_{I,exact}$ in a calculus context and $\tau_{P,prim}$ in a probability context and whether they know technologies related to IC and CPD. Then, by questions 6 and 7, we focus on relationships between CPD and IC and the previous questions: we analyse links between the questions 3, 4 (second part), 5 and 6, on the one hand, and links between the questions 2 and 7, on the other hand.

**DATA ANALYSIS AND RESULTS**

The questionnaire was used at the beginning of September 2017 (the first week of classes) in two classes of first year CPGE (French engineers school preparatory classes) students, which is in fact at the transition between secondary and tertiary levels. The first class (called class N) is a class of MPSI (Mathematics, Physics and Engineering Science) and the second class (called class R) is a class of PCSI (Physics, Chemistry and Engineering Science) of a rather prestigious establishment. The students working on the questionnaire are in selective classes, so we can assume that they are "good" scientific students, and in particular, if they meet difficulties then these are shared by the other students. We only analysed answers from students who studied in French high school during the past year because we constructed the questionnaire taking into account the context of the French high school institution. We retrieved 82 questionnaires (40 of the class N and 42 of the class R). Except for a few students (less than 5), the students didn't use a calculator during the test.

From the 82 students, only one does not mobilize the technique $\tau_{I,exact}$ to resolve the routine task concerning IC (question 1). 85% find a correct expression of the primitive and 78% obtain the correct result, which means that this technique is mastered quite well by students. Regarding the computation of the probability for an exponential distribution (question 4), 63% of the students get a correct result. 82% of the students use the technique $\tau_{P,prim}$ and the other students directly apply a formula. More than 70% of the students could identify the probability as an integral. Summarizing, except for some errors regarding the primitive or the computation, the majority of students are able to pass from a probability to an integral and, moreover, know the fundamental theorem of calculus. In praxeological terms, they are able to mobilize the technique $\tau_{P,prim}$ that implies the use of the technique $\tau_{I,exact}$ in the case of an exponential distribution.

Regarding CPD and neglecting formulation and formalization issues, only 39% of the students know the definition of the density function (question 2) and only 27% recollect the density function of the normal distribution (question 3b). In view of question 5, only 32% mobilize the theorem claiming that all continuous functions on an interval admit a primitive. 52 % give an example of a function for which they don't explicitly know a primitive (although it might exist and be expressed in terms of standard functions[2], for instance $ln(x)$). Among these, 15% mention the density function of the Gaussian distribution (or a function of the type $e^{-x^2}$). Of the 22 students who know the density function of the normal distribution, half propose it as an example in question 5 (13% of all the students).

9% of the students do not answer question 6. The method most often proposed is technique $\tau_{I,exact}$ (46%). Less than 20% (16 students) mention the normal distribution and 16% propose the rectangle method ($\tau_{I,approx}$). A few students propose the technique "integration by parts", which is beyond the curriculum in grade 12. Moreover, only 23% of the students who propose the rectangle method are able to illustrate the method by drawing the graph of the Gaussian curve (the other students draw a wrong curve or do not consider any graph). Only 11 of the 16 students (69%) who mention the normal distribution in this question write that the integral $I$ is equal to the probability $P(-0,5 \leq X \leq 1)$ with $X$ a random variable of a reduced normal distribution and 7 of them state more precisely that they have to use the calculator to evaluate this probability ($\tau_{P,norm}$). So, our results indicate that the ostensive $\int_{-0,5}^{1} \frac{e^{-x^2/2}}{\sqrt{2\pi}} dx$ without indication invite students to stay in the IC sector and even more particularly in the praxeology $\Pi_{I,exact}$ (even if it is not possible here), which reflects that $\Pi_{I,exact}$ is the praxeology most developed in grade 12. To pass from an integral to a probability and to change the sector in this direction does not seem to be natural for students. Moreover, of the 22 students who know the expression for the density function of the reduced normal distribution (question 3b), only 10 (45%) recognize the density function of the normal distribution is this context. Of the 15 students who say that the Gaussian function does not admit primitives expressed by using standard elementary functions (question 5b), around 27% proposes to use the technique $\tau_{I,exact}$ nevertheless and only 53% recognize the normal distribution in question 6. This means that most of the students were not able to mobilize $\Pi_{P,norm}$ in the IC-context of question 6, that is $\Pi_{I,norm}$. In particular the application of $\tau_{I,exact}$, although question 5b is answered correctly, demonstrates the strong compartmentalisation between CPD and IC.

Regarding the answers to question 7, we notice that 38% of the students (31) state that the function $f$ is a density function. 42% of them (13 students) could justify that it is a density function including 4 who forget to mention the positivity of the function (because they don't write this condition in their definition in question 2). Moreover, 12 of the 31 students as well as one additional one identify $A(\lambda)$ as $P(0 \leq X \leq \lambda)$. 8 students think that the distribution in question 7 is an exponential distribution of parameter $x$ and 7 students state that $A$ is a density function. Finally, more than 46% of the students mention "probability", which means that they manage to identify at least some link between the IC embedding of question 7 and CPD. Probably, the questionnaire itself influenced students and the percentage would be lower otherwise, i.e. if question 7 was asked independently. Overall correct results with justification are rare and of the 32 students who master the definition of a probability density function (one element of the technology $\theta_{P,prim}$ tested in question 2), 31% (10 students) do not identify $f$ as a density function. 16% of the students do not at all answer question 7.

Again, we observe also with respect to question 7 that CPD praxeologies, although they are in principle available, often cannot be mobilized in the IC contextualization.

Summarizing, the data analysis shows that techniques related to one sector are available for the majority of students only when ostensives related to this sector are provided. Additionally, related technologies are much less mastered by students. Perhaps, this could be an explanation why it is not natural for students to mobilize praxeologies of the CPD sector for a task in the IC sector, in addition to the fact that these tasks are not taught in the classroom. The data analysis by all means shows a strong compartmentalisation between CPD and IC.

## CONCLUSION AND PERSPECTIVES

The results of our primary data analysis clearly demonstrate a strong compartmentalisation between CPD and IC. In particular, techniques from CPD, although available in a CPD task, could not be mobilized in an IC-contextualized task. A next step in our research will be more detailed data analyses looking for correlations and interdependencies between techniques and technologies of CPD and IC. We further observed that available techniques were not accompanied by related technologies. One could claim that more elaborated technologies might support the transfer of techniques. More generally, we think about studies investigating the impact of changes in the institutional setting, i.e. establishing innovative teaching sequences with less compartmentalisation. A teaching sequence articulating CPD and IC is proposed in Derouet (2016). The effect of this teaching on the answers of students to the questionnaire could be analysed by comparing the latter with the present results.

## NOTES

1. This axiomatic remains implicit at the secondary level; it may be related to Measure Theory at university level.

2. We realised *a posteriori* that our question was not phrased properly: "a function that admits a primitive whose expression is not explicitly known to you" may be interpreted as a lack of techniques to actually compute the primitive and not the impossibility to provide an expression (in terms of elementary functions).